\documentclass[12pt]{amsart}
\usepackage{amscd,amsthm,amssymb,amsfonts,amsmath,epsfig,epic,bbm,url}
\usepackage[arrow,matrix,tips,2cell]{xy}

\theoremstyle{plain}
\newtheorem{theorem}{Theorem}
\newtheorem{lemma}[theorem]{Lemma}

\theoremstyle{definition}
\newtheorem{definition}[theorem]{Definition}

\theoremstyle{remark}
\newtheorem*{remark}{Remark}

\newcommand{\Z}{\mathbb Z}    % Integers
\newcommand{\R}{\mathbb R}    % Reals
\newcommand{\C}{\mathbb C}    % Complexes
\newcommand{\Aff}{\operatorname{Af{}f}}

\newcommand{\Pic}{\operatorname{Pic}}

\newcommand{\Div}{\operatorname{Div}}

\newcommand{\length}{\operatorname{length}}

\newcommand{\val}{\operatorname{val}}

\newcommand{\dist}{\operatorname{dist}}

\newcommand{\ignore}[1]{\relax}
\newcommand{\dd}{\partial}

\begin{document}

\title{Tropical Theta characteristics}

\author{Ilia Zharkov}
\address{Kansas State University, 138 Cardwell Hall, Manhattan, KS 66506}
\email{zharkov@math.ksu.edu}

\begin{abstract}
This note is a follow up of  \cite{MZ2} and it is largely inspired by a beautiful description of \cite{BN} of non-effective degree $g-1$ divisors via chip-firing game. We consider the set of all theta characteristics on a tropical curve and identify the Riemann constant $\kappa$ as a unique non-effective one among them. 
\end{abstract}
\maketitle

First we recall the basic setup and main players. For more details the reader can look in \cite{MZ2} and references there in. 

Let $\Gamma$ be a connected finite graph and ${\mathcal V}_1(\Gamma)$ be the set of its 1-valent vertices. We say $\Gamma$ is a metric graph if the topological space $\Gamma\setminus{\mathcal V_1}(\Gamma)$ is given a complete metric structure and $\Gamma$ is its compactification. In particular, all leaves have infinite lengths. Introducing a new interior two-valent vertex on any edge is set to give an equivalent metric graph, and a {\em tropical curve $C$} is an equivalence class of such graphs. Its genus is $g=b_1(\Gamma)$ for any representative $\Gamma$.

The metric allows one to talk about affine and piece-wise linear functions on $C$ with integral slopes. At every vertex $v$ we may define the set of outward primitive tangent vectors $\xi_i$. Then any PL function $f$ defines a {\em principal divisor} 
$$(f)=\sum_{p\in C} (\sum_{i=1}^{\val(p)}\frac{\dd f(p)}{\dd \xi_i}) p. 
$$
In general, a {\em divisor} $D=\sum a_i p_i$ is a formal linear combination of points in $C$ with integral coefficients. We say: $D_1\sim D_2$ if $D_1-D_2$ is a principal divisor, $D\geq 0$ if all $a_i\geq 0$, and $D$ has degree $\deg D=\sum a_i$. The degree $2g-2$ divisor $K:=\sum_{p\in C} (\val(p)-2)p$ is called the {\em canonical divisor}.

Let $\Aff$ be the sheaf of $\Z$-affine functions (i.e., local PL functions $f$ with $(f)=0$). Define the integral cotangent local system
${\mathcal T_\Z}^*$ on $C$ by the following exact sequence of sheaves:
\begin{equation*}
0 \longrightarrow \R  \longrightarrow \Aff \longrightarrow {\mathcal T_\Z}^*\longrightarrow
0,
\end{equation*}
The rank $g$ lattice of {\em 1-forms} $\Omega_\Z(C)$ (a.k.a. the {\em circuit lattice}) on $C$ is formed by the global sections of ${\mathcal T_\Z}^*$. Let $\Omega(C)^*$ be the vector space of $\R$-valued linear functionals on $\Omega_\Z(C)$. Then the integral cycles $H_1(C,\Z)$ form a lattice $\Lambda$ in $\Omega(C)^*$ by integrating over them. We define the tropical {\em Jacobian} (cf. \cite{Nagnibeda2}) to be
$$J(C):= \Omega(C)^*/H_1(C,\Z)\cong \R^g/\Lambda.$$

Let us fix a reference point $p_0\in C$. Given a divisor
$D=\sum a_i p_i$ we choose paths from $p_0$ to $p_i$. Integration along these
paths defines a linear functional on $\Omega_\Z(C)$:
$$\hat\mu(D)(\omega)=\sum a_i \int_{p_0}^{p_i} \omega. $$
For another choice of paths the value of $\hat\mu(D)$ will differ by an element
in $\Lambda$. Thus, we get a well-defined tropical analog of the {\em
Abel-Jacobi map} $\mu^d:\Div^d(C)\to J(C)$. The tropical Abel-Jacobi theorem (cf. \cite{Nagnibeda2},\cite{BN},\cite{MZ2}) says that
for each degree $d$ the map $\mu^d$ factors through $\Pic^d(C)$ (the group of divisors modulo linear equivalence):
\begin{center}
\mbox{}\xymatrix{ \Div^d (C)\ar[dr]_{\mu^d} \ar[r] & \Pic^d (C) \ar[d]^{\phi}\\
     & J(C)}
\end{center}
and $\phi$ is a bijection. Both maps $\mu^d$ and $\phi$ depend on the base point $p_0$ unless $d=0$.

An explicit solution to the Jacobi inversion problem is provided by introducing the tropical theta function as follows. The metric on $C$ defines a symmetric positive bilinear form $Q$ on $\Omega(C)^*$ by setting $Q(\ell,\ell) := \length(\ell)$ on simple cycles $\ell$. That, in turn, defines a convex $\Lambda$-quasi-periodic PL function on $\Omega(C)^*$:
\begin{equation}\label{eq:theta}
\Theta(x):=\max_{\lambda\in\Lambda} \{Q(\lambda,x) - \frac12
Q(\lambda,\lambda)\}, \quad x\in\Omega(C)^*,
\end{equation}
that can be thought as a section of a polarization line bundle on $J(C)$. Its corner locus defines the {\em theta divisor} $[\Theta]$ on the Jacobian.

For $\lambda\in J(C)$ let $[\Theta_\lambda]=[\Theta]+\lambda$ be the translated theta divisor on $J(C)$  and let $D_\lambda:=\mu^*[\Theta_\lambda]$ denote the (effective, of degree $g$) pull back divisor of
$[\Theta_\lambda]$ to the curve via the Abel-Jacobi map $\mu:C\to J(C)$.

\begin{theorem}[Jacobi inversion, \cite{MZ2}]\label{jacobi_inversion}
There exists a universal $\kappa\in \Pic^{g-1}(C)$ such that $\mu^g(D_\lambda)+\kappa=\lambda$ for all $\lambda\in J(C)$.
\end{theorem}

Let $W_{g-1} \subset \Pic^{g-1}(C)$ denote the Abel-Jacobi image of the set of effective divisors of degree $g-1$. 

\begin{theorem}[Riemann's theta divisor \cite{MZ2}]\label{thm:riemann}
$W_{g-1}+\kappa=[\Theta]$, and $2(-\kappa)=\mu^{2g-2}(K)$.
\end{theorem}

\begin{definition}
Divisor classes $\mathcal K \in\Pic^{g-1}(C)$ such that $2 \mathcal K =\mu^{2g-2}(K)$ are called theta characteristics. We set $\mathcal K_0=-\kappa$.
\end{definition}

\begin{lemma}
Among the set of $2^g$ theta characteristics $\mathcal K_0+\frac12\Lambda/\Lambda\subset \Pic^{g-1}(C)$ only $\mathcal K_0$ is not in $W_{g-1}$.
\end{lemma}
\begin{proof}
According to the Theorem \ref{thm:riemann} it suffices to prove that among the two-torsion points $\frac12\Lambda/\Lambda\subset J(C)$ only $0$ is not in $[\Theta]$. But since the form $Q$ is not degenerate the maximum at $x=0$ in (\ref{eq:theta}) is achieved by a single term, namely by $\lambda=0$. Hence $0\in\Omega(C)^*$ is not in the corner locus of $\Theta$.

On the other hand if a non-zero $\lambda\in \frac12\Lambda$ belongs to the interior of the maximal Voronoi cell containing 0, then so does $-\lambda$ (because $\Theta(x)$ is an even function). But this is impossible since $\lambda-(-\lambda)=2\lambda\in \Lambda\setminus\{0\}$. 
\end{proof}

\begin{remark}
The fact that one theta characteristic is special in tropical world stems from choosing a Lagrangian splitting of the lattice  $H_1(C_\C,\Z)$ when degenerating classical curves. Definition of the classical theta function depends on this choice. Once the choice is made there is a special even element among the classical theta characteristics (cf. e.g. \cite{Mumford}). 
\end{remark}

Next we will give an explicit geometric description of the theta characteristics. Let $S\subset C$ be a finite subset. Consider the distance function $d_S(x)=\dist(S,x)$. Its gradient flow defines an {\em acyclic} orientation on (some graph $\Gamma$ in the class of) $C$. Let 
$$\mathcal K^{\pm}_S= \sum_{p\in C} (\val_{\pm}(p)-1)p
$$
be the two associated moderators (cf. \cite{BN}, \cite{MZ2}). Here $\val_{\pm}$ stands for the number of outgoing/incoming edges.

\begin{lemma}
$\mu^{g-1}(\mathcal K^{\pm}_S)=-\mathcal K_0$.
\end{lemma} 
\begin{proof}
Both $\mathcal K^{\pm}_S$ are divisors of degree $g-1$ which represent non-effective classes in $\Pic^{g-1}(C)$ (cf. \cite{BN}, \cite{MZ2}). Because $\mathcal K^+_S+\mathcal K^-_S=K$ it suffices to show that $\mathcal K^+_S\sim \mathcal K^-_S$. In fact, we claim $\mathcal K^+_S=\mathcal K^-_S+(d_S)$.

Indeed, at any $p\in C$ the slope of $d_S$ along an outward primitive vector $\xi_i$ is +1 if $\xi_i$ agrees with the gradient flow of $d_S$, and -1 if it does not. Thus, $\deg(d_S)\vert_p=\val_+(p)-\val_-(p)$, which put together gives 
$\mathcal K^+_S-\mathcal K^-_S=  \sum_{p\in C} (\val_+ -\val_-)p = (d_S).$
\end{proof}

As a consequence we note that for all subsets $S$ the corresponding $\mathcal K^{\pm}_S$ represent the same class $\mathcal K_0$. Explicitly, $\mathcal K^-_{S'}-\mathcal K^-_S=\frac12(d_{S'}-d_{S})$. The minimal representatives (having a single pole) of $\mathcal K_0$ are given by $S=q$, a point in $C$. Then $\mathcal K^-_q$ is the $q$-reduced form of $\mathcal K_0$ (cf. \cite{BN}, \cite{MZ2}). 

Next we identify effective theta characteristics.  For each non-trivial element $\gamma\in H_1(C,\Z/2\Z)=\frac12\Lambda/\Lambda$ we choose a simple representative $\bar \gamma \in H_1(C,\Z)$, that is, a cycle containing any edge of $C$ with multiplicity at most one. Let $|\gamma|$ be the support of $\bar\gamma$. It depends only on $\gamma$ (not on $\bar\gamma$), contains only finite edges and has only even-valent vertices. Let $d_\gamma(x):=\dist(|\gamma|,x)$. The gradient flow of $d_\gamma$ together with the cyclic orientation of $\bar\gamma$ give an orientation on $C$, which, in turn, defines a pair of divisors
$$\mathcal K^{\pm}_\gamma= \sum_{p\in C} (\val_{\pm}(p)-1)p.$$
Again, $\mathcal K^{\pm}_\gamma$ depends only on $\gamma$, not on its lift $\bar\gamma$. Since the inducing orientation is not acyclic, the $\mathcal K^{\pm}_\gamma$ are not moderators. In fact, $\mathcal K^{-}_\gamma$ is manifestly effective because the flow has no source point.

Similar to the $\mathcal K^{\pm}_S$ we have $\mathcal K^{+}_\gamma+\mathcal K^{-}_\gamma=K$ and $\mathcal K^{+}_\gamma - \mathcal K^{-}_\gamma = (d_\gamma)$. Thus every $\mathcal K^{\pm}_\gamma$ represents a theta characteristic. We will write $\mathcal K_\gamma$ for the class of $\mathcal K^{\pm}_\gamma$ in $\Pic^{g-1}(C)$.

\begin{lemma}
$\mathcal K_\gamma-\mathcal K_0 = \frac12\gamma$ in $J(C)$.
\end{lemma}
\begin{proof}
First, note that $\frac12\gamma$ is well defined in $J(C)$, i.e. independent of the choice of $\bar\gamma$. Next we take $S$ to be any finite subset of $|\gamma|$ containing all boundary points of $C\setminus|\gamma|$ and all vertices of $|\gamma|$. Then $|\gamma|\setminus S$ consists of disjoint open intervals, and we denote by $M$ the set of their mid points.

Clearly $d_\gamma=d_S$ on $C\setminus|\gamma|$, hence $\mathcal K^{\pm}_\gamma=\mathcal K^{\pm}_S$ on $C\setminus|\gamma|$.  On the other hand, on $|\gamma|$ we have 
$$
\mathcal K^{-}_\gamma = \sum_{p\in S} (\frac 12 \val_\gamma (p)-1) p, \quad \text{and} \quad 
\mathcal K^-_S = M - S,
$$
where $\val_\gamma(p)$ denotes the valence of $p$ as a vertex of $|\gamma|$. Put together, 
$$
\mathcal K^-_\gamma-\mathcal K^-_S= \sum_{p\in S} ( \frac 12 \val_\gamma (p) ) p-M.
$$ 
Now any lift $\bar\gamma$ of $\gamma$ specifies a system of paths by starting at points of $M$ and ending at points of $S$. The linear functional on $\Omega(C)$ defined by these paths is precisely equal to $\frac12\bar\gamma$, and its projection to $J(C)$ is $\frac12\gamma$. On the other hand,  
since there are exactly $\frac 12 \val_\gamma (p)$ paths ending at every vertex $p\in S$ this functional  represents the divisor $\mathcal K^-_\gamma-\mathcal K^-_S$.
\end{proof}

We conclude by combining all results in a single statement.
\begin{theorem}
The map $\gamma\mapsto \mathcal K_\gamma$ from $H_1(C,\Z/2\Z)$ to the set of theta characteristics is an isomorphism of affine $\Z/2\Z$-lattices. Among all $\mathcal K_\gamma$ only $\mathcal K_0=-\kappa$ is not (linearly equivalent to) an effective divisor.
\end{theorem}

\end{document}